\documentclass[12pt]{article}
\newif\ifhide
\hidefalse

\usepackage{amsmath}
\usepackage{amssymb}
\usepackage{color}
\newtheorem{Theorem}{Theorem}
\newtheorem{Lemma}{Lemma}
\newtheorem{Corollary}{Corollary}
\newtheorem{Remark}{Remark}
\newtheorem{Proposition}{Proposition}
\newtheorem{Definition}{Definition}

\title{{\normalsize\tt\hfill\jobname.tex}\\
Note on 
local mixing techniques for stochastic differential equations
 \\
}
\author{A. Yu. Veretennikov\footnote{School of Mathematics,
University of Leeds, LS2 9JT, UK, 
HSE National university, Moscow,  and Institute of
Information Transmission Problems, Russian Academy of
Sciences, Moscow, Russian Federation,  E-mail:  {\tt\small
a.veretennikov @ leeds.ac.uk}}}

\begin{document}

\maketitle

%Author's commands
%\newtheorem{theorem}{Theorem}

%\date{}
%\newtheorem{Example}{Example}
%\newtheorem{Theorem}{Theorem}
%\newtheorem{Lemma}{Lemma}
%\newtheorem{Coro\exp(-ct)llary}{Corollary}
%\newtheorem{Remark} {Remark}
%\newtheorem{Proposition}{Proposition}
%g in a step-wise fashion,
%uncomment
% the following command:

%\beamerdefaultoverlayspec

%\begin{document}

 % \titlepage
 
\begin{abstract}
This paper discusses several techniques which may be used for applying the coupling method to solutions of stochastic differential equations (SDEs). They all work in dimension $d\ge 1$, although, in $d=1$ the most natural way is to use intersections of trajectories, which requires nothing but strong Markov property and non-degeneracy of the diffusion coefficient.   
In dimensions $d>1$ it is possible to use embedded Markov chains either by considering discrete times $n=0,1,\ldots$, or by arranging special stopping time sequences and to use 
local Markov -- Dobrushin's (MD) condition. Further applications may be based on one or another version of the MD condition. For  studies of  convergence and mixing rates  the (Markov) process must be strong Markov and recurrent; however, recurrence is a separate issue which is not discussed in this paper.
\end{abstract}

%\section{Markov--Dobrushin's condition via Girsanov}

  \section{Introduction}
%\begin{abstract}
The stochastic differential equation (SDE) in $R^d$
\begin{equation}\label{sde1}
X_t^{} = x + \int_{0}^t b(X^{}_s)ds + \int_{0}^t \sigma(X^{}_s)dW_s, \, t\ge 0, 
\end{equation}
is considered. Here $(W_t), t\ge 0$ is a $d$-dimensional Wiener process, $b$ and $\sigma$ are vector and matrix-valued Borel measurable functions of dimensions $d$ and $d\times d$ respectively.

It is assumed that the equation (\ref{sde1}) has a (weak or strong) solution which is a strong Markov process; see \cite{Krylov_selection}. Naturally, under this condition the process $X_n$ -- that is, our solution $X_t$ considered at integer times $t=0,1,\ldots $ -- is a Markov chain (MC), which is, of course, also strong Markov. The advantage of the total variation distance (although, it is not unique in this respect) for Markov processes is that once it is established that, for example, 
$$
\| \mu_n - \mu \|_{TV} \le \psi(n) \to 0, \quad n\to \infty, 
$$
where $\mu_t$ is the marginal distribution of $X_t$, $\mu$ is any probability measure (ergodic limit for $(\mu_n)$),   
then this rate of convergence can be nearly verbatim transported to the continuous time:
$$
\| \mu_t - \mu \|_{TV} \le \psi([t]) \to 0, \quad t\to \infty, 
$$
where $[t]$ is the integer part of $t$. 
Consider two independent versions of our Markov process $X_t$ (in continuous time), say, $X_t^2$ and $X_t^2$, with two different initial values $x^1$ and $x^2$, respectively (or distributions). Since we allow weak solutions, the processes $X_t^2$ and $X_t^2$, generally speaking, are defined on two different probability spaces with two different Wiener processes; without loss of generality, we may assume that they are {\em independent} and consider their direct product; thus, we have two trajectories $X^2$ and $X^2$ on the same probability space (do not forget that Wiener processes are also different and independent of each other). Denote by $Q(x,dx')$ the transition kernel, 
$$
Q(x,dx') = \mathbb P_x(X_1 \in dx').
$$
\begin{Definition}
A (global) Markov -- Dobrushin's (MD) condition holds for the Markov process $X_n$ iff 
\begin{equation}\label{MD1}
\inf_{x_1,x_2} \int \left(\frac{Q(x_1,dx')}{Q(x_2,dx')}\wedge 1\right)Q(x_2,dx') > 0.
\end{equation}
\end{Definition}
Here it is not assumed that necessarily $Q(x_1,dx')$ is absolutely continuous with respect to $Q(x_2,dx')$, but its absolutely continuous component and its derivative is taken under the integral. In what follows a localised version of this condition will be stated and this localised version will be the object of our main interest in this paper.
General approaches to coupling for SDEs require a (usually positive) recurrence and some form of local mixing. For the latter, beside intersections applicable only in the case $d=1$, the following tools may be used. 
\begin{itemize}
\item
Lower and upper transition density bounds (requires H\"older coefficients and ``elliptic'' or ``hypoelliptic'' non-degeneracy); here the popular in discrete time theory of ergodic Markov chains {\em petite sets condition} along with recurrence properties {\em may} be used.

\item
Lower and upper density bounds of the  transition density {\em only for the equation without the drift,} including degenerate and highly degenerate cases: here Girsanov's transformation is an efficient tool; petite sets conditions, generally speaking, do not work.

\item
In the absence of lower and upper density bounds, under the non-degeneracy condition and for general measurable coefficients of SDE Harnack inequalities in parabolic or elliptic versions may be applied; a petite sets condition apparently may be proved; however, they are less efficienet than MD because the latter guarantees better convergence rate esimates. 

\end{itemize}
Hence, the goal of this paper is to attract more attention to the MD condition, which this condition deserves in the humble opinion of the author. There is also a hope that this list of available techniques may help in the future in studying ergodic properties for more general classes of processes. 

One more point is that except for the method based on lower and upper density bounds, in all other more involved situations the popular DD condition \cite{Doob}, or, in its local version, the {\em petite set} condition is difficult to apply to SDEs, unlike the MD one; and even if it may be applied, the MD condition provides weaker assumptions and better convergence rate estimates. 
%From the formal point of view, the petite set condition is a uniform minorization of all transition kernels from any state in some suitable set to the same set by a unique minorization measure. At the same time, MD condition is designed for comparison of any two such kernels with each other directly, without any auxiliary minorization measure. Hence, the latter condition (MD) is more general. 

Note that for discrete time stochastic models, and in certain cases for continuous time, too, one more natural approach to coupling is to use regeneration. Unfortunately, for general SDEs this method is not available. So, we do not discuss it here, although, the multidimensional coupling constructions for processes with continuous distributions are sometimes called ``a generalised regeneration''.

Let us warn the reader that most of the results of this paper are known, perhaps, in a slightly different form; we just collect them here together.  For simplicity we do not touch more general equations such as SDEs with jumps. However, in principle, more general Markov processes, in particular, SDEs with jumps, may also be tackled with the help of similar techniques.  
  
It should be also highlighted that all methods discussed in what follows (but the elliptic Harnack inequality) can be applied with minor differences to non-homogeneous SDEs, too, except that convergence would be for the distance in total variation between marginal distributions corresponding to any two initial measures, not to the invariant measure which does not exist in this situation, in general.

The paper consists of two sections: this introduction and the main section 2; in turn section 2 is split into six sub-sections, most of them related to one of  the coupling tools listed above. The majority of  proofs are sketchy or dropped because the results are known; the only exception is the part about elliptic Harnack, which is new to the best of the author's knowledge and where the details of the proof are shown. The paper presents the set of various {\em tools} for coupling for SDEs. Neither recurrence -- the necessary second ingredient in studying convergence and mixing rates -- nor coupling itself (except for the basic lemma \ref{odvuh} added for the reader's convenience) are not the goals of this paper.

\section{Main results}

\subsection{Case $d=1$, MD \& coupling using intersections}

If in the 1D case for local coupling we can use {\bf intersections} of two independent solutions of the same SDE with different initial values. Assume that $X_t$ and $X'_t$ are two solutions of the equation (\ref{sde1}) with different initial values $X_0=x$ and $X'_0=x'$ in the one-dimensional case. The basis for applying coupling via intersections is the following result. 

\begin{Proposition}\label{Pro1}
If $b$, $\sigma$, and $\sigma^{-1}$ are bounded then 
\begin{equation}
\inf_{-1\le x,x'\le 1}P_{x,x'}(\exists \, s\in[0,1]: \, X_s=X'_s)>0.  
\end{equation}
The first meeting time $\tau:= (t\ge 0: X_t=X'_t)$ is a stopping time.
\end{Proposition}
{\em Proof} follows from the following two elementary steps. 

1. Change time for both SDEs making diffusion coefficients equal to one. There is no need to make it {\em the same} random time change: generally speaking, the latter is not possible unless the diffusion coefficient is a constant.
% in the case of the variable diffusion coefficient.) 
Since $\sigma$ and $\sigma^{-1}$ are bounded, the interval $[0,1]$ after this change becomes random, but for both processes contains some non-random interval $[0,T]:=[0,\inf_x \sigma^{-2}(x)]$. This can be also applied to non-homogeneous SDEs with coefficients depending on time. 

2. The random time change leaves the drift bounded. Hence, due to Girsanov's transformation of measure, the probability that the process with a lower initial value will attain the level $+1$ over $[0,T]$ is positive and bounded away from zero. 
Similarly, the probability that the process with a bigger initial value will attain the level $-1$ over $[0,T]$ is positive and bounded away from zero. Therefore, they meet on $[0,T]$ with a positive probability which is bounded away from zero, as required. \hfill QED

\subsection{MD, ``case b'' \& ``petite set'' conditions}
In dimensions $d>1$ intersections do not work for the ``normal'' SDEs, and we now switch to the main topic of this paper, local mixing conditions. 
Global and local versions of ``petite set'' and Markov--Dobrushin's (MD) conditions will be stated. Most frequently either of them is applied in its local variant, but the global options also work in cases of a uniform ergodicity. It should be noted that, in fact, local versions may vary slightly depending on a particular setting; we only show their main appearances. The ``petite set'' condition is a localised version of the ``case (b)'' condition from \cite[Chapter V, section 5]{Doob}, which is a simplification of the ``condition D'' (nowadays called Doeblin -- Doob's one) from the same chapters in \cite{Doob}. Let us highlight that the MD condition may also be in a global or local form.

\begin{Definition}
The process satisfies the condition ``b'' (from \cite[Chapter V]{Doob}) iff there exists a probability measure $\nu$ on the state space $\cal X$ and constants $T,c>0$ such that 
\begin{equation}\label{grandeset}
%\inf_{x_0, x_1} \int_{\cal X} \left(\frac{\mu^{x_0}_T(dy)}{\nu(dy)}\wedge \frac{\mu^{x_1}_T(dy)}{\nu(dy)}\right){\nu(dy)} >0.
\inf_{x_0\in \cal X}
\frac{\mu^{x_0}_T(dy)}{\nu(dy)}\ge c.
\end{equation}
\end{Definition}

\begin{Definition}
The process satisfies the local condition `b', or ``petite set'' condition iff there exists a set $D \subset {\cal X}$, a  probability measure $\nu$ on $D$ and constants $T,c>0$ such that 
\begin{equation}\label{grandeset}
%\inf_{x_0, x_1} \int_{\cal X} \left(\frac{\mu^{x_0}_T(dy)}{\nu(dy)}\wedge \frac{\mu^{x_1}_T(dy)}{\nu(dy)}\right){\nu(dy)} >0.
\inf_{x_0\in D}
\frac{\mu^{x_0}_T(dy)}{\nu(dy)}\ge c.
\end{equation}
\end{Definition}
See, in particular, \cite{MeynTweedie} about the usage of the petite set condition in convergence studies. 
Recall that normally this local condition -- as well as the local MD condition in the definition \ref{def4} in the next paragraphs -- should be accomplished with certain recurrence assumptions or properties; however, as it was said earlier, recurrence is not the goal of this paper, and it makes sense to study it separately. Global conditions ``case b'' and MD both lead to efficient uniform in the initial data  exponential convergence.

\begin{Definition}
The following is called the global Markov--Dobrushin condition: there exists $T>0$ such that 
\begin{equation}\label{kappaglobal}
\kappa(T):=\inf_{x_0, x_1} 
\int_{} \left(\frac{\mu^{x_0}_T(dy)}{\mu^{x_1}_T(dy)}\wedge 1\right)\mu^{x_1}_T(dy) >0.
\end{equation}

\end{Definition}

\begin{Definition}\label{def4}
The following is called a local Markov--Dobrushin condition: there exist sets $D, D'\subset \cal X$ in the state space and a constant $T>0$ such that 
\begin{equation}\label{kappalocal}
\kappa(D,D';T):=\inf_{x_0, x_1 \in D} 
\int_{D'} \left(\frac{\mu^{x_0}_T(dy)}{\mu^{x_1}_T(dy)}\wedge 1\right)\mu^{x_1}_T(dy) >0,
\end{equation}
\end{Definition}

\begin{Remark}
Usually, but not necessarily $D'=D$; in this case we use the notation $\kappa(D,D';T) =: \kappa(D;T)$. Another possibility is $D' = R^d$. 
A sufficient condition for (\ref{kappalocal}) is as follows: there exists a dominating measure $\nu(dy)$ such that $\mu^{x}_T(dy) \ll \nu(dy)$ for any $x \in D$, and 
\begin{equation}\label{kappalocal2}
\kappa(D,D';T)=\inf_{x_0, x_1 \in D} 
\int_{D'} \left(\frac{\mu^{x_0}_T(dy)}{\nu(dy)}\wedge \frac{\mu^{x_1}_T(dy)}{\nu(dy)}\right)\nu(dy) >0.
\end{equation}
In general, there might be no dominating measure for all $x$ simultaneously. Yet, as we shall see, (\ref{kappalocal2}) may be verified in most of the cases in what follows.

\end{Remark}
Clearly, the ``petite set'' condition implies the MD one, both in the global (``case b'') and local versions, for example, (\ref{grandeset}) implies (\ref{kappaglobal}):
\begin{align*}%\label{kappaglobal}
\inf_{x_0, x_1} 
\int_{} \left(\frac{\mu^{x_0}_T(dy)}{\mu^{x_1}_T(dy)}\wedge 1\right)\mu^{x_1}_T(dy) 
 %\\\\
= \int_{} \left(\frac{\mu^{x_0}_T(dy)}{\nu(dy)}\wedge \frac{\mu^{x_1}_T(dy)}{\nu(dy)}\right)\nu(dy) \ge c > 0, 
\end{align*}
but, generally speaking, not vice versa. The basis for applying coupling via any of  them is the following {\em coupling lemma} (not to be confused with the {\em coupling inequality}). Let us add that the MD condition admits some further generalisation, see \cite{VV, VV2}, which provides in certain cases a slightly better efficient  convergence bound under slightly wider assumptions. However, this  note is just about tools which allow to check  a local condition MD for non-degenerate SDEs. The following lemma clarifies why MD condition is so useful; at the same time it serves as the basis for a further application of the MD condition to coupling technique for Markov processes. 

\begin{Lemma}[``Of two random variables'']\label{odvuh}
Let $X^{1}$ and $X^2$ be two random variables on their (without loss of generality different, which will be made independent after we take their direct product) probability spaces $(\Omega^1, {\cal F}^1, \mathbb P^1)$ and $(\Omega^2, {\cal F}^2, \mathbb P^2)$ and with densities $p^1$ and $p^2$ with respect to some reference measure $\Lambda$, correspondingly.  Then, if 
\begin{equation*}%\label{MDkappa}
1-p:= q = \int \left(p^1(x)\wedge p^2(x)\right) \Lambda(dx) > 0, 
\end{equation*}
then there exists one more probability space $(\Omega, {\cal F}, \mathbb P)$ and two random variables on it $\tilde X^1, \tilde X^2$ such that 
\begin{equation*}%\label{eq2rv}
\!\!{\cal L}(\tilde X^j)\! =\!{\cal L}(X^j), j\!=\!1,\!2, \; \& \;  \frac{\|{\cal L}(X^1) \!-\! {\cal L}(X^2)\|_{TV}}2 \! = \!P(\tilde X^1 \!\not=\! \tilde X^2) \!=\! p. 
\end{equation*}
 \end{Lemma}
This is a well-known technical tool in the coupling method. The proof -- which is simple enough -- may be found, for example, in \cite{Ver_ln}. This reference should not be regarded as a  claim that this lemma belongs to the author, although, who whe first inventor of this lemma is not clear to him. 
%Formally, it is probably due to \cite{Griffeath}, or \cite{Nummelin}; informally, its idea may belong to Doob \cite{Doob}, or to Doeblin \cite{Doeblin}?, or to  Kolmogorov \cite{Kolm_MP}, or even to Markov himself \cite{Markov}. 

The next lemma justifies the hint that to estimate the convergence rate  for a Markov process to its invariant measure (assume that it exists) in continuous time $(X_t, t\ge 0)$ it suffices to evaluate it for discrete times $n=0,1,\ldots$ Its elementary proof is provided for the reader's convenience. 
Let $\mu^X_t$ be the marginal distribution of $X_t$, and let  $\mu^X$ be its invariant measure.  
\begin{Lemma}\label{lem2}
It holds, 
$$
\|\mu^X_t - \mu^X\|_{TV} \le \|\mu^X_n - \mu^X\|_{TV}, \quad t\ge n.
$$

\end{Lemma}
{\em Proof.} Due to the Markov's property of $X$, by Chapman -- Kolmogorov's equation  using the convention $a_+ = a\vee 0$, $a_- = (-a)\vee 0$, we get 
\begin{align*}
\frac12\|\mu^X_t - \mu^X\|_{TV} = \sup_A
(P_x(X_t\in A) - P_\mu(X_t\in A))
 \\
= \sup_A\iint 1(z\in A) (P_x(X_n\in dy) - P_\mu(X_0\in dy))P_y(X_{t-n}\in dz)
 \\
= \sup_A\left(\iint 1(z\in A) (P_x(X_n\in dy) - P_\mu(X_0\in dy))_+P_y(X_{t-n}\in dz) \right.
 \\
\left. - \iint 1(z\in A) (P_x(X_n\in dy) - P_\mu(X_0\in dy))_-P_y(X_{t-n}\in dz)\right) 
 \\
\le \sup_A \iint 1(z\in A) (P_x(X_n\in dy) - P_\mu(X_0\in dy))_+P_y(X_{t-n}\in dz)
 \\
= \iint P_y(X_{t-n}\in dz) (P_x(X_n\in dy) - P_\mu(X_0\in dy))_+
 \\
= \int  (P_x(X_n\in dy) - P_\mu(X_0\in dy))_+
= \frac12\|\mu^X_n - \mu^X\|_{TV},
\end{align*}
%and
%\begin{align*}
%\iint 1(z\in A) (P_x(X_n\in dy) - P_\mu(X_t\in dy))_-P_y(X_{t-n}\in dz)\le \frac12\|\mu^X_n - \mu^X\|_{TV}.
%\end{align*}
as required. \hfill QED

\ifhide
 
{\em Proof.}
 {Assume $0< q<1$, otherwise the lemma is trivial.}
Let r.v. $\eta_1$, $\eta_2$, $\xi$, have the following densities with respect to the measure $\Lambda$:
\begin{align*}
p_{\eta^1}(t) &= (1-q)^{-1}\left(p^1(t) - p^1(t) \wedge
 p^2(t)\right),
\\
p_{\eta^2}(t)&= (1-q)^{-1}\left(p^2(t) - p^1(t) \wedge
 p^2(t)\right),
\\
p_{\xi}(t) &= q^{-1}\left(p^1(t) \wedge p^2(t)\right).
\end{align*}
Let $\zeta$ be a random variable independent of $\eta^1$,
 $\eta^2$ and $\xi$ taking values in $\{0,1\}$ such
that
\begin{equation*}
P(\zeta=0)=q, \;\; P(\zeta=1)=1-q.
\end{equation*}
Assume that $q\neq0$ and $q\neq1$ and let 
\begin{align*}
\tilde X^1:=\eta^1 1(\zeta=1)+\xi 1(\zeta=0),
\\
\tilde X^2:=\eta^2 1(\zeta=1)+\xi 1(\zeta=0).
\end{align*}
Then $\tilde X^1\stackrel{d}{=}X^1$, $\tilde
 X^2\stackrel{d}{=}X^2$, and $P(\tilde
X^1=\tilde X^2)= q$, \quad QED.

\fi

\subsection{MD using lower and upper transition density bounds}
Assume $d>1$ and Gaussian type upper and lower bounds for the transition densities, which can be established under the non-degeneracy of $\sigma\sigma^*$ and H\"older coefficients \cite{Eidelman, Friedman, Solo}, as well as under certain hypoellipticity conditions (see, for example, \cite{DelarueMenozzi, Menozzi, Polidoro}, et al.). {\em In particular,}  let $\sigma\sigma^*$ be uniformly non-degenerate, and let both $\sigma$ and $b$ satisfy H\"older's conditions: there exists $L, \alpha>0$ such that for any $x,x'$
\begin{align*}\label{hoelder}
|b(x)-b(x')| + \|\sigma(x)-\sigma(x')\| \le L |x-x'|^\alpha. 
\end{align*}
As it follows from the PDE theory (see the references above), under such conditions for any $t>0$ there exist constants $C_t,C'_t,c_t,c'_t>0$ such that Gaussian type lower and upper bounds hold true for the transition densities $f_t(x,x')$ (fundamental solutions in the PDE language)
\begin{align*}%\label{hoelder}
C'_t\exp(-c_t |x-x'|^2)\le f_t(x,x') \le C_t \exp(-c^{-1}_t |x-x'|^2).
\end{align*}
In particular, under the non-degeneracy condition on $\sigma\sigma^*$, it may be used $C_t = C t^{-d/2}$,  $C'_t = C' t^{-d/2}$, $c_t = c t$, $c'_t = c' t$ with some $C, C', c, c'$, and under the hypoelliptic conditions it is also known how to evaluate all these constants. 
Then, clearly, a local ``petite set'' condition is satisfied with any bounded domain $D$ (an open set by definition) and with the Lebesgue measure as $\nu$. 
Hence, the MD condition is also valid. 
To the best of the author's knowledge this is the only case -- although, this class of coefficients is wide enough, but far from the most general -- where the ``petite set'' condition can be applied to Markov SDEs in order to arrange coupling.

\medskip

%For simplicity, we restrict our study to the case $\sigma \equiv I$, a unit diffusion matrix. In this case (as well as in the non-degenerate case!) it turns out that we may use Girsanov's transformation instead of Harnack. An alternative simplified way would be to use PDE results on the transition densities  satisfying Gaussian type lower and upper bounds [Solonnikov, Eidelman, Friedman]; it will be also mentioned briefly in the end. 

\subsection{MD using stochastic exponentials}
In this section let us assume that lower and upper bounds for transition densities hold true for the SDE with a ``truncated drift''  
$$
X^0_t = x + \int_{0}^t \sigma(X^{0}_s)dW_s +\int_{0}^t b_1(X^{0}_s)ds,
$$
while the goal is to arrange local coupling for the full SDE with the more involved drift of the form 
$$b= b_1 + b_2,$$
where $b_2$ is just Borel measurable and bounded (this boundedness may be relaxed). 
We are interested in establishing an MD condition for the full equation  (\ref{sde1}). Note that upper and lower bounds from the previous subsection, in general, are not applicable. 
Denote 
$$
\tilde b_2(x) := \sigma^{-1}(x) b_2(x),
$$
and let
$$
\rho_T \!:= \!\exp\!\left(\! -\! \int_0^T \!\tilde b_2(X_t) \, d W_t
 \!-\!\frac12 \!\int_0^T\!
\left|\tilde b_2(X_t)\right|^2 \, dt \right)\!.
$$
Recall that $\rho_T$ is a probability density for any $T>0$. Denote by $\mu_t$ the marginal distribution of $X_t$.
 \begin{Theorem}[local MD condition via Girsanov]\label{thm1}
For any $T>0$ and $R>0$
\begin{equation}\label{kappa}
\kappa(R,T):=\inf_{x_0, x_1\in B_R} 
\int_{B_R} \left(\frac{\mu^{x_0}_T(dy)}{\mu^{x_1}_T(dy)}\wedge 1\right)\mu^{x_1}_T(dy) >0.
\end{equation}
 \end{Theorem}
This inequality suffices for applications to coupling and  convergence rates (given suitable recurrence estimates). For the {\em proof} of very close statements (actually, even for degenerate SDEs) see \cite{Abourashchi, Ver_AiT20}. Some other localised versions of this result may be established: as an example, the sets $B_R$ under the infimum sign and as a domain of integration may, actually, differ.

\subsection{MD using parabolic Harnack inequalities}
As usual in this paper, in this section we assume that $d\ge 1$, coefficients $b$ and $\sigma$ are bounded (which can be relaxed by a localisation) and Borel measurable, and  $\sigma\sigma^*$ is uniformly non-degenerate.
Under such conditions Krylov -- Safonov's Harnack parabolic inequality holds true \cite[Theorem 1.1]{KS_parab}, stated here in terms of probabilities rather than solutions of PDEs:
\begin{align}\label{ks_parab}
\sup_{|x_1|, |x_2|\le 1/4}\frac{P_{}(X^{0,x_1}_\tau \in d\gamma)}{P_{}(X^{\epsilon,x_2}_\tau \in d\gamma)}|_{\Gamma_\epsilon}\le N<\infty, 
\end{align}
where $\Gamma_\epsilon$ is the parabolic boundary of the cylinder $((t,x): |x|\le 1; \epsilon \le t\le 1)$, i.e. ($\Gamma_\epsilon = \Gamma^{(t=1)}_\epsilon \cup \Gamma^{(t<1)}_\epsilon$),
$\Gamma_\epsilon = ((t,x): (|x|=1 \; \& \; \epsilon \le t\le 1) \cup (|x|\le 1 \; \& \; t=1))$, and 
$$
\tau:= \inf(t\ge 0: |X_t|\ge 1), \; \text{with a convention}\; \inf(\emptyset)=1;
$$
the constant $N$ depends on $d$, on the ellipticity constants of the diffusion, on the sup-norm of the drift, and on $\epsilon$. Note that in (\ref{ks_parab}) the measure in the numerator is absolutely continuous with respect to the one in the denominator, that is, there is no singular component in this situation. 
Let 
$$
\mu^{x_1}(d\gamma) = P_{}(X^{0,x_1}_\tau \in d\gamma),
\quad \mu^{\epsilon, x_2}(d\gamma) = P_{}(X^{\epsilon,x_2}_\tau \in d\gamma), 
$$ 
where $d\gamma$ is the element of the boundary $\Gamma_\epsilon$. 
Then the following local mixing bound holds true.
\begin{Theorem}[local MD via parabolic Harnack]\label{thm2}
Let $\mu^{x_1}(\Gamma_\epsilon)\ge q$ with some $q>0$. Then a version of Markov-Dobrushin's condition holds, 
\begin{equation}\label{MDparab}
\inf_{|x_1|, |x_2|\le 1/4} 
\int_{\Gamma_\epsilon} \left(\frac{\mu^{\epsilon, x_2}(d\gamma)}{\mu^{x_1}(d\gamma)}\wedge 1\right)\mu^{x_1}(d\gamma) \ge \frac{q}{N} >0.
\end{equation}
\end{Theorem}
Note that the value $q$ here may be chosen arbitrarily close to one, if $\epsilon>0$ is small enough. However, the decrease of  $\epsilon$ implies the increase of the constant $N$ in (\ref{ks_parab}).

\medskip

\noindent
{\bf Proof}.  
Indeed, due to the inequality (\ref{ks_parab}) we have, 
$$
f := \frac{d\mu^{x_1}}{d\mu^{\epsilon,x_2}}|_{\Gamma_\epsilon}\le N \quad \& \quad
\mu^{x_1}\ll  \mu^{\epsilon, x_2} \; \text{on $\Gamma_\epsilon$}  
$$
Denote by $\tilde \mu^{\epsilon,x_2}$ the absolutely continuous part of $\mu^{\epsilon,x_2}$ with respect to $\mu^{x_1}$ (we do not know whether there exists a singular component here, but the calculus in what follows does not depend on this). 
%on $\Gamma_\epsilon$. 
Then 
$$
\frac{d\tilde \mu^{\epsilon,x_2}}{d\mu^{x_1}}|_{\Gamma_\epsilon} = \frac1{f}\ge \frac1{N}.
$$
Hence, the assumption $\mu^{x_1}(\Gamma_\epsilon) \ge q$ implies 
\begin{align*}
\int_{\Gamma_\epsilon} \left(\frac{\mu^{\epsilon,x_2}(d\gamma)}{\mu^{x_1}(d\gamma)}\wedge 1\right)\mu^{x_{x_1}}(d\gamma) 
=\int_{\Gamma^{}_\epsilon} \left(\frac{\tilde \mu^{\epsilon,x_2}(d\gamma)}{\mu^{x_1}(d\gamma)}\wedge 1\right)\mu^{x_1}(d\gamma) 
 \\\\
\ge \int_{\Gamma^{}_\epsilon} \frac1{N} \mu^{x_1}(d\gamma)
\ge \frac{q}{N} >0. \hspace{3cm} \text{QED}
\end{align*}
%Note that there is no hint for establishing the petite set condition under the assumptions of the theorem here, due to the difference in the definitions of the measures $\mu_1$ and $\mu_2$, as the second one also depends on some $\epsilon>0$: so, they are not the transition kernels over the same interval of time from two different states. Nevertheless, this version of the MD condition does allow an efficient construction of coupling. Technically, it happens because the MD condition (as well as the petite set condition) requires only a one-side inequality, not bounds from above and from below.

%\medskip

Sometimes it may be more convenient to use another version of the MD condition, which follows from theorem \ref{thm2}. Denote 
$$
\mu^x_{1}(A):= P_x(X_{1} \in A), 
%\quad \Lambda_z^{}(A) :=P_z(X_{1-\epsilon} \in A),   
\quad
 A\subset R^d.
$$

\begin{Corollary}\label{cor1}
Under the assumptions of theorem \ref{thm2} 
the following version of the MD condition holds: there exists $q'\in (0,q)$ such that
\begin{align}\label{newellMD}
\inf_{|x_1|,|x_2|\le 1/8}\int_{R^d} \left(\frac{\mu^{x_1}_{1}(dy)}{ \mu^{x_2}_{1}(dy)}\wedge 1\right) \mu^{x_2}_{1}(dy) \ge 
\frac{q'}{N}.
\end{align}
\end{Corollary}
Note that here $R^d$ plays the role of $D'$ in the MD condition. In some cases this may not be convenient; however, using moment bounds of the solution a reasonable version of this inequality with some bounded ball $B_R$ in place of $R^d$ is, of course, possible. We leave it till further studies where such a replacement may be required. 

\medskip

\noindent
{\bf Proof.} Note that due to the boundedness of $\sigma$ and $b$, 
$$
\inf_{|x|\le 1/8}P_x(|X_\epsilon|\le 1/4)> 0.
$$
Denote
$$
q' := q \inf_{|x|\le 1/8}P_x(|X_\epsilon|\le 1/4).
$$
We have, 
\begin{align*}
\mu^{x_2}_{1}(dy)=P_{x_2}(X_1\in dy) 
= E_{x_2} E(X_1\in dy |X_\epsilon) = 
 \\\\ 
\ge E_{x_2} 1(|X_\epsilon| \le 1/4)
 E(X_1\in dy |X_\epsilon) = E_{x_2} 1(|X_\epsilon| \le 1/4)
 \mu_1^{\epsilon, X_\epsilon}(dy).
\end{align*}
Hence, denoting $\nu_{\epsilon,x_2}(dz):= P_{x_2}(X_\epsilon \in dz)$, we find
\begin{align*}
1 \wedge \frac{ \mu^{x_2}_{1}(dy)}{\mu^{x_1}_{1}(dy)} 
 \ge 1 \wedge \frac{E_{x_2} 1(|X_\epsilon| \le 1/4)
 \mu_1^{\epsilon, X_\epsilon}(dy)}{\mu^{x_1}_{1}(dy)}
 \\\\
=  1 \wedge \frac{\int \nu_{\epsilon, x_2}(dz)1(|z| \le 1/4)
 \mu_1^{\epsilon, z}(dy)}{\mu^{x_1}_{1}(dy)} 
 \\\\
= 1 \wedge \left( \int \nu_{\epsilon, x_2}(dz)1(|z| \le 1/4)\frac{
 \mu_1^{\epsilon, z}(dy)}{\mu^{x_1}_{1}(dy)} \right)
 \\\\
\ge  \int \nu_{\epsilon, x_2}(dz)\left(1 \wedge 1(|z| \le 1/4) \frac{ \mu_1^{\epsilon, z}(dy)}{\mu^{x_1}_{1}(dy)}\right)
 \\\\
\ge  \int \nu_{\epsilon, x_2}(dz)1(|z| \le 1/4) \left(1 \wedge\frac{ \mu_1^{\epsilon, z}(dy)}{\mu^{x_1}_{1}(dy)}\right).
\end{align*}
So, 
\begin{align}\label{mdchar}
\int \left(1 \wedge \frac{ \mu^{x_2}_{1}(dy)}{\mu^{x_1}_{1}(dy)}\right) \mu^{x_1}_{1}(dy)
  \nonumber \\\nonumber\\ \nonumber
\ge \int \left[\int \nu_{\epsilon, x_2}(dz)1(|z| \le 1/4) \left(1 \wedge\frac{ \mu_1^{\epsilon, z}(dy)}{\mu^{x_1}_{1}(dy)}\right)\right] \mu^{x_1}_{1}(dy) 
  \\ \nonumber \\ 
=\int \nu_{\epsilon, x_2}(dz)1(|z| \le 1/4) \left[\int \left(1 \wedge\frac{ \mu_1^{\epsilon, z}(dy)}{\mu^{x_1}_{1}(dy)}\right) \mu^{x_1}_{1}(dy)\right].
\end{align}
This was the first step in the reduction of the MD characteristics in the left hand side of (\ref{mdchar}) to the Harnack inequality: now we may deal with the measures   $\mu_1^{\epsilon, z}(dy)$ and $\mu^{x_1}_{1}(dy)$. However, these are still not the ones which show up  in (\ref{ks_parab}) or in (\ref{MDparab}). The next step will complete this reduction. Let 
$$
\tilde \Lambda_{x_1,z}(dy) := \mu_1^{\epsilon, z}(dy) + \mu^{x_1}_{1}(dy). 
$$
Then 
\begin{equation}\label{e14}
\int \left(1 \wedge\frac{ \mu_1^{\epsilon, z}(dy)}{\mu^{x_1}_{1}(dy)}\right) \mu^{x_1}_{1}(dy) = \int \left(\frac{ \mu_1^{x_1}(dy)}{\tilde \Lambda_{x_1,z}(dy)} \wedge \frac{ \mu_1^{\epsilon, z}(dy)}{\tilde \Lambda_{x_1,z}(dy)}\right) \tilde \Lambda_{x_1,z}(dy).
\end{equation}
We have, 
\begin{align}\label{e15}
\!\!\!\!\!\frac{\mu^{\epsilon,z}_{1}(dy)}{\tilde \Lambda_{x_1,z}(dy)}\wedge \frac{ \mu^{x_1}_{1}(dy)}{\tilde \Lambda_{x_1,z}(dy)} 
\! \ge \! 
\frac{P_{x_1}(X_1\in dy, \tau < 1)}{\tilde \Lambda_{x_1,z}(dy)}\wedge \frac{P_{z}(X_{1-\epsilon}\in dy, \tau < 1-\epsilon)}{\tilde \Lambda_{x_1,z}(dy)}. 
\end{align}
Therefore, 
\begin{align*}%\label{newellMD}
\int \left(\frac{ \mu_1^{x_1}(dy)}{\tilde \Lambda_{x_1,z}(dy)} \wedge \frac{ \mu_1^{\epsilon, z}(dy)}{\tilde \Lambda_{x_1,z}(dy)}\right) \tilde \Lambda_{x_1,z}(dy)
 \\\\
\ge  \int_{R^d} \left(\frac{P_{x_1}(X_1\in dy, \tau < 1)}{\tilde \Lambda_{x_1,z}(dy)}\wedge \frac{P_{z}(X_{1-\epsilon}\in dy, \tau < 1-\epsilon)}{\tilde \Lambda_{x_1,z}(dy)}\right)\tilde \Lambda_{x_1,z}(dy)
 \\\\
= \int_{R^d} \left(\frac{P_{x_1}(X_1\in dy, \tau < 1)}{P_{z}(X_{1-\epsilon}\in dy, \tau < 1-\epsilon)}\wedge 1\right) P_{z}(X_{1-\epsilon}\in dy, \tau < 1-\epsilon).
\end{align*}
Further, since $|x_1|\le 1/4$ and $|z|\le 1/8$, then 
%by Chapman -- Kolmogorov's equation, and 
due to the strong Markov propery and by virtue of the inequality (\ref{ks_parab}) we have, 
%for any $0<t<1-\epsilon$, 
\begin{align*}
P_{x_1}(X_1\in dy, \tau < 1)  
= E_{x_1} 1( \tau < 1) (X_1\in dy) 
 \\\\
=E_{x_1} E\left(1( \tau < 1) (X_1\in dy)|{\cal F}_{\tau}\right) 
=E_{x_1} 1( \tau < 1) E\left((X_1\in dy)|{\cal F}_{\tau}\right) 
 \\\\
= E_{x_1} 1( \tau < 1) E\left((X_1\in dy)|X_{\tau}\right) 
 \\\\
= E_{x_1} 1( \tau < 1) E_{t,y}\left((X_{1-t}\in dy)\right)|_{(t,y)=(\tau,X_{\tau})} 
 \\\\
\ge \frac{q}{N} E_{x_2} 1( \tau < 1) E_{t,y}\left((X_{1-t}\in dy)\right)|_{(t,y)=(\tau,X_{\tau})} 
 \\\\
= \frac{q}{N}  P_{z}(X_{1-\epsilon - t}\in dy, \tau < 1).
\end{align*}
So, 
\begin{align*}
\int_{R^d} \left(\frac{P_{x_1}(X_1\in dy, \tau < 1)}{P_{z}(X_{1-\epsilon}\in dy, \tau < 1)}\wedge 1\right) P_{z}(X_{1-\epsilon}\in dy, \tau < 1)
 \\\\
\ge \frac{q}{N} \int_{R^d} P_{z}(X_{1-\epsilon}\in dy, \tau < 1) 
= \frac{q}{N} P_{z}(\tau < 1).
\end{align*}
Recall that in (\ref{mdchar}) the integrand involves the indicator $1(|z|\le 1/4)$. Clearly, 
$$
\kappa:= \inf_{|z|\le 1/4} P_{z}(\tau < 1) >0. 
$$
Hence, due to (\ref{mdchar}), (\ref{e14}) and (\ref{e15}), 
\begin{align*}%\label{e16}
\int \nu_{\epsilon, x_2}(dz)1(|z| \le 1/4) \left[\int \left(1 \wedge\frac{ \mu_1^{\epsilon, z}(dy)}{\mu^{x_1}_{1}(dy)}\right) \mu^{x_1}_{1}(dy)\right] 
 \\\\
\ge \frac{q\kappa}{N} \int \nu_{\epsilon, x_2}(dz)1(|z| \le 1/4) 
=  \frac{q\kappa}{N} \int 1(|z| \le 1/4) P_{x_2}(X_\epsilon \in dz) \ge \frac{q'}{N}
\end{align*}
with some $0<q'<q\kappa$, 
as required. The inequality (\ref{newellMD}) follows.  \hfill QED

\subsection{MD using elliptic Harnack inequalities}
The assumptions of this sections are the same as in the previous one: $d\ge 1$, coefficients $b$ and $\sigma$ are bounded (which can be relaxed) and Borel measurable, and  $\sigma\sigma^*$ is uniformly non-degenerate. We have the elliptic Harnack inequality due to \cite[Theorem 3.1]{Safonov}, stated here in its probabilistic form (while in \cite{Safonov} it is offered in the language of elliptic PDEs): there exists a constant $N>0$ such that for any $0< R \le 1$ and any $A \in \partial B_R$,
\begin{equation}\label{Ha-ell}
\sup_{|x|\le R/8} \mathbb P_x(X_{\tau_R} \in A) 
\le N \inf_{|x|\le R/8} \mathbb P_x(X_{\tau_R} \in A), 
\end{equation}
where $\tau_R = \inf(t\ge 0: |X_t|\ge R)$, and $\partial B_R$ is the boundary of the ball $B_R$. This inequality itself is some MD condition. In fact, it is not clear whether this version of Harnack inequality may be helpful for estimating convergence rate of the distribution of $X_t$ to its stationary regime.   Nevertheless, if it can be used for such a purpose -- which is the author's hope -- then it might be more convenient to apply  the following version of the MD condition based on the inequality (\ref{Ha-ell}). 
Let
$$
Q_{R,T} := \{(t,x): \; t\le T, \; |x|\le R\}.
$$
Note that 
$$
\bigcup_{T > 0}Q_{R,T} = R_+ \times B_R = R_+ \times (x: \, |x|\le R).
$$
Denote by $\Gamma_{R,T}$ the part of the parabolic boundary of $Q_{R,T}$ corresponding to $t<T$, namely, 
$$
\Gamma_{R,T} = ((t,x): \, 0\le t\le T, \, |x|=R).
$$
Denote $\tau_{R,T}:= \inf(t\ge 0: X_t\not\in B_{R}) \wedge T$, 
$\tau_{R}:= \inf(t\ge 0: X_t\not\in B_{R})$, and 
$$
\nu^x_{R,T}(A):= P_x(X_{\tau_{R,T}} \in A), \quad
\nu^x_R(A):= P_x(X_{\tau_{R}} \in A), \quad A\subset \partial B_R.
$$

\begin{Theorem}[local MD via elliptic Harnack]\label{thm3}
The local MD condition  
\begin{align}\label{MDT0}
\inf_{x_1,x_2\in B_R}\int_{B_R} \left(\frac{ \nu^{x_1}_{R,T}(d\gamma)}{\nu^0_R (d\gamma)}\wedge \frac{ \nu^{x_2}_{R,T}(d\gamma)}{\nu^0_R (d\gamma)}\right)\nu^0_R (d\gamma) \ge 
\frac{C_R}{2N}
\end{align}
holds true for any $T>0$ large enough.
\end{Theorem}
In principle, it is possible to evaluate such values of $T$ for which (\ref{MDT0}) holds, and it might be useful for estimating convergence rates, but, as it was said earlier,  we do not pursue this goal here.    
As already mentioned, unlike with the parabolic Harnack inequality, it is not clear how useful the local mixing property in the form (\ref{MDT0}) of (\ref{Ha-ell}) could be for coupling; this may be clarified in further studies. 

\medskip

\noindent
{\bf Proof.} 
Clearly, $\tau_{R,T} \le \tau_{R}$. Note  that due to the non-degeneracy of $\sigma$ we have $\tau_R<\infty$ a.s., and 
$$
%\lim_{T} P(|X_{\tau_{R,T}}|<R) = 0, \;\&\; 
\lim_{T\to\infty} \inf_{|x|\le R/8}P_x(\tau_{R,T}=\tau_R) 
= \lim_{T\to\infty} \inf_{|x|\le R/8}P_x(\tau_{R}<T) = 1.
$$
Equivalently, 
$$
%\lim_{T} P(|X_{\tau_{R,T}}|<R) = 0, \;\&\; 
\lim_{T\to\infty} \sup_{|x|\le R/8}P_x(\tau_{R,T}<\tau_R) 
= \lim_{T\to\infty} \inf_{|x|\le R/8}P_x(T<\tau_{R}) = 0.
$$
Hence, 
\begin{align}\label{Tuni}
\nu^x_{R,T}(A)= P_x(X_{\tau_{R,T}} \in A) = P_x(X_{\tau_{R}} \in A, \tau_R < T) \uparrow 
P_x(X_{\tau_{R}} \in A) = \nu^x_R(A), 
\end{align}
as $T\uparrow \infty$, where the convergence is uniform with respect to $A$ and $|x|\le R/8$.

The inequality (\ref{Ha-ell}) implies the following, 
\begin{align*}
0< N^{-1 }\le \inf_{|x|\le 1/8}\frac{ \nu^x_R(d\gamma)}{\nu^0_R (d\gamma)}|_{\partial B_R}\le \sup_{|x|\le 1/8}\frac{ \nu^x_R(d\gamma)}{\nu^0_R (d\gamma)}|_{\partial B_R}\le N<\infty.
\end{align*}
%or, equivalently, for any $A\subset \partial B_R$
%\begin{align*}
%0< N^{-1 }\le \inf_{|x|\le 1/8}\frac{ \nu^x_R(A)}{\nu^0_R (A)}|_{\partial B_R}\le \sup_{|x|\le 1/8}\frac{ \nu^x_R(A)}{\nu^0_R (A)}|_{\partial B_R}\le N<\infty, 
%\end{align*}
As a consequence, for any $R>0$ we get
\begin{align}\label{CRN}
\inf_{x_1,x_2\in B_R} \int_{B_R} \left(\frac{ \nu^{x_1}_R(d\gamma)}{\nu^0_R (d\gamma)}\wedge \frac{ \nu^{x_2}_R(d\gamma)}{\nu^0_R (d\gamma)}\right)\nu^0_R (d\gamma) \ge \frac{C_R}{N}
>0.
\end{align}
By virtue of the monotone convergence theorem and due to (\ref{Tuni}) 
we have for any $x_1,x_2$, 
\begin{align*}
%\inf_{x_1,x_2\in B_R} 
\int_{B_R} \left(\frac{ \nu^{x_1}_{R,T}(d\gamma)}{\nu^0_R (d\gamma)}\wedge \frac{ \nu^{x_2}_{R,T}(d\gamma)}{\nu^0_R (d\gamma)}\right)\nu^0_R (d\gamma) 
 \\\\ 
\to 
%\inf_{x_1,x_2\in B_R} 
\int_{B_R} \left(\frac{ \nu^{x_1}_R(d\gamma)}{\nu^0_R (d\gamma)}\wedge \frac{ \nu^{x_2}_R(d\gamma)}{\nu^0_R (d\gamma)}\right)\nu^0_R (d\gamma), \quad T\to\infty.
\end{align*}
Hence, for $T$ large enough we obtain from (\ref{CRN}),

%for any $A\subset \partial B_R$, 
%\begin{align*}%\label{HRT}
%0< \frac12 N^{-1 }\le \inf_{|x|\le 1/8}\frac{\nu^{x}_{R,T}(A)}{\nu^0_R (A)}|_{\partial B_R}\le \sup_{|x|\le 1/8}\frac{\nu^{x}_{R,T}(A)}{\nu^{0}_{R}(A)}|_{\partial B_R}\le 2N<\infty, 
%\end{align*}
%or, equivalently, 
\begin{align*}
%\inf_{x_1,x_2\in B_R} 
\int_{B_R} \left(\frac{ \nu^{x_1}_{R,T}(d\gamma)}{\nu^0_R (d\gamma)}\wedge \frac{ \nu^{x_2}_{R,T}(d\gamma)}{\nu^0_R (d\gamma)}\right)\nu^0_R (d\gamma) 
\ge \frac{C_R}{2N}. 
\end{align*}
%as required. \hfill QED
However, technically this is still not sufficient because we want a similar inequality with infimum $\inf_{x_1,x_2\in B_R}$. Using the elementary inequality $(a-b)\wedge (c-d) \ge a\wedge c - b - d$ along with the identity
$$
 P_x(X_{\tau_{R}} \in d\gamma, \tau_R < T) =  P_x(X_{\tau_{R}} \in d\gamma) -  P_x(X_{\tau_{R}} \in d\gamma, \tau_R \ge T),
$$
we have,

\begin{align*}
%\inf_{x_1,x_2\in B_R} 
\int_{B_R} \left(\frac{ \nu^{x_1}_{R,T}(d\gamma)}{\nu^0_R (d\gamma)}\wedge \frac{ \nu^{x_2}_{R,T}(d\gamma)}{\nu^0_R (d\gamma)}\right)\nu^0_R (d\gamma) 
 \\\\
=\int_{B_R} \left(\frac{ P_{x_1}(X_{\tau_{R}} \in d\gamma, \tau_R < T)}{\nu^0_R (d\gamma)}\wedge \frac{P_{x_1}(X_{\tau_{R}} \in d\gamma, \tau_R < T)}{\nu^0_R (d\gamma)}\right)\nu^0_R (d\gamma) 
 \\\\
\ge  \int_{B_R} \left(\frac{ P_{x_1}(X_{\tau_{R}} \in d\gamma)}{\nu^0_R (d\gamma)}\wedge \frac{P_{x_2}(X_{\tau_{R}} \in d\gamma)}{\nu^0_R (d\gamma)}\right)\nu^0_R (d\gamma)
 \\\\
- \int_{B_R} \left(\frac{ P_{x_1}(X_{\tau_{R}} \in d\gamma, \tau_R \ge T)}{\nu^0_R (d\gamma)} + \frac{P_{x_2}(X_{\tau_{R}} \in d\gamma, \tau_R \ge  T)}{\nu^0_R (d\gamma)}\right)\nu^0_R (d\gamma).
\end{align*}
Here 
\begin{align*}
( \sup_{x_1,x_2\in B_R})\; \int_{B_R} \left(\frac{ P_{x_1}(X_{\tau_{R}} \in d\gamma, \tau_R \ge T)}{\nu^0_R (d\gamma)} + \frac{P_{x_2}(X_{\tau_{R}} \in d\gamma, \tau_R \ge  T)}{\nu^0_R (d\gamma)}\right)\nu^0_R (d\gamma)
 \\\\
=( \sup_{x_1,x_2\in B_R})\; \int_{B_R} \left(P_{x_1}(X_{\tau_{R}} \in d\gamma, \tau_R \ge T) + P_{x_2}(X_{\tau_{R}} \in d\gamma, \tau_R \ge  T)\right) 
 \\\\
\le 2 \sup_{x\in B_R} P_{x}(\tau_R \ge T)\to 0, \quad T\to\infty.  
\end{align*}
On the other hand, 
\begin{align*}
\inf_{x_1,x_2\in B_R}\int_{B_R} \left(\frac{ P_{x_1}(X_{\tau_{R}} \in d\gamma)}{\nu^0_R (d\gamma)}\wedge \frac{P_{x_2}(X_{\tau_{R}} \in d\gamma)}{\nu^0_R (d\gamma)}\right)\nu^0_R (d\gamma) \ge \frac{C_R}{N}.
\end{align*}
So, 
\begin{align*}
\inf_{x_1,x_2\in B_R}\int_{B_R} \left(\frac{ \nu^{x_1}_{R,T}(d\gamma)}{\nu^0_R (d\gamma)}\wedge \frac{ \nu^{x_2}_{R,T}(d\gamma)}{\nu^0_R (d\gamma)}\right)\nu^0_R (d\gamma) 
 \\\\ 
\to 
\inf_{x_1,x_2\in B_R} 
\int_{B_R} \left(\frac{ \nu^{x_1}_R(d\gamma)}{\nu^0_R (d\gamma)}\wedge \frac{ \nu^{x_2}_R(d\gamma)}{\nu^0_R (d\gamma)}\right)\nu^0_R (d\gamma), \quad T\to\infty.
\end{align*}
Therefore, there exists $T_0$ such that 
%for any $T\ge T_0$ we have 
\begin{align}\label{MDT0}
\inf_{x_1,x_2\in B_R}\int_{B_R} \left(\frac{ \nu^{x_1}_{R,T}(d\gamma)}{\nu^0_R (d\gamma)}\wedge \frac{ \nu^{x_2}_{R,T}(d\gamma)}{\nu^0_R (d\gamma)}\right)\nu^0_R (d\gamma) \ge 
\frac{C_R}{2N}, \quad T\ge T_0,
\end{align}
which completes the proof. \hfill QED

\section*{Acknowledgements}
The article  was prepared within the framework of the HSE University Basic Research Program in part which includes Theorem \ref{thm1}  and all lemmata. The part of Theorems \ref{thm2} and \ref{thm3} and Corollary \ref{cor1} was funded by Russian Science Foundation grant 17-11-01098.
 
 %{\em \mbox{https://rdcu.be/b4h3F}} and in the preprint 
%{\em \mbox{https://arxiv.org/abs/2006.12134}}

%\verbatim {https://rdcu.be/b4h3F}

\end{document}